\documentclass{article}

\usepackage{graphicx}
\usepackage{multirow}
\usepackage{listings}
\usepackage{amsmath}
\usepackage{amsfonts}
\usepackage{caption}
\usepackage[margin=1in]{geometry}
\usepackage{soul,color}
\usepackage{cite}
\usepackage[caption=false]{subfig}

\definecolor{dblue}{RGB}{33, 0, 75}

\begin{document}

\markboth{V. S. Duzhin, N. N. Vasilyev}
{Asymptotic behaviour of normalized dimensions of standard and strict Young diagrams \,--\, growth and oscillations}

\title{Asymptotic behaviour of normalized dimensions of standard and strict Young diagrams \,--\, growth and oscillations}

\author{V. Duzhin, N. Vasilyev}
\date{}

\maketitle

    \parbox{.95\textwidth}{%
\begin{footnotesize} \textcolor{dblue} {\textbf{For citation}: V.~S.~Duzhin, N.~N.~Vasilyev, \emph{Asymptotic behavior of normalized dimensions of standard and strict Young diagrams --- growth and oscillations.} --- J. Knot Theory Ramifications \textbf{26} (2016). doi:10.1142/S0218216516420025} \end{footnotesize}
     }%

\begin{abstract}
In this paper, we present the results of a computer investigation of asymptotics for maximum dimensions of linear and projective representations of the symmetric group. This problem reduces to the investigation of standard and strict Young diagrams of maximum dimensions. We constructed some sequences for both standard and strict Young diagrams with extremely large dimensions. The conjecture that the limit of normalized dimensions exists was proposed by A. M. Vershik 30 years ago \cite{verker85} and has not been proved yet. 

We studied the growth and oscillations of the normalized dimension function in sequences of Young diagrams.
Our approach is based on analyzing finite differences of their normalized dimensions.
This analysis also allows us to give much more precise estimation of the limit constants.
\end{abstract}

\hfill \begin{minipage}{0.275\textwidth}
\emph{
Dedicated to memory of\\Sergei Duzhin (1956 -- 2015).
}
\end{minipage}

\section{Introduction}

\subsection{Definitions}

\begin{figure}[h]
	\centering

	\subfloat[][]{\includegraphics[scale=0.6]{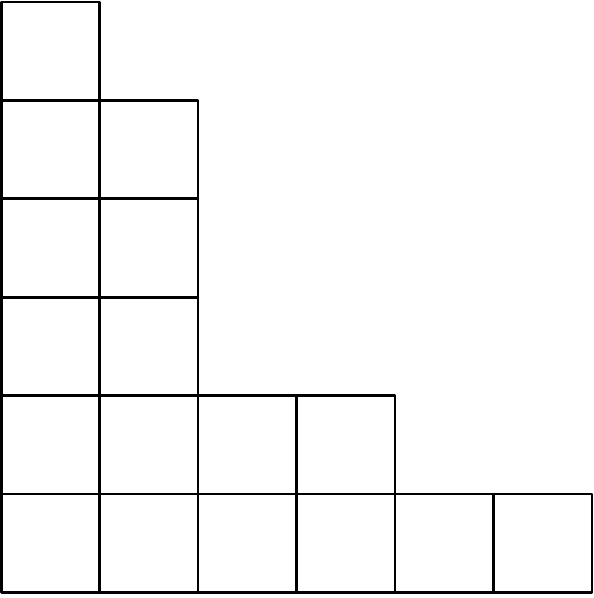}\label{fig:standard}}
\hspace{100pt}
     	\subfloat[][]{\includegraphics[scale=0.6]{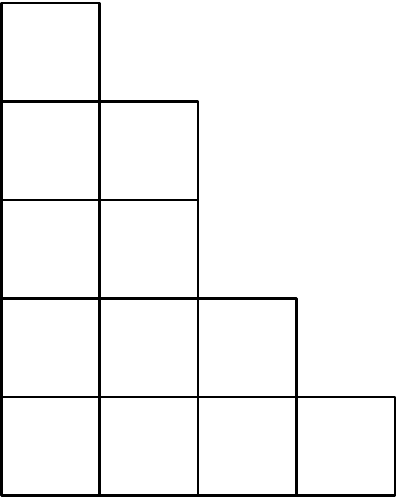}\label{fig:strict}}
     	\caption{Examples of Young diagrams: (a) a Standard Young diagram, (b) a strict Young diagram.}
     	\label{fig:2diags}

	\end{figure}

A Young diagram of the shape $\lambda=(l_1,l_2,\ldots,l_k)$ consists of $k$ columns. The column numbered as $i$ consists of $l_i$ boxes. There are several equivalent ways to present Young diagrams. In this paper, we draw the diagrams in the upper right quadrant. Figure \ref{fig:2diags} illustrates two examples of Young diagrams represented by the partitions (6, 5, 2, 2, 1, 1), (5, 4, 2, 1), respectively.

A \textit{strict Young diagram} is a Young diagram with no columns of the same height. Figure~\ref{fig:standard} is standard, while Fig.~\ref{fig:strict} is strict. The strict diagrams are related to so-called strict partitions of natural numbers. A partition $({l_1,l_2,\ldots,l_k})$ of a number $n$ is called \textit{strict} if it does not contain components of the same size.

\begin{figure}[!htb]
	\centering
	\includegraphics[scale=0.3]{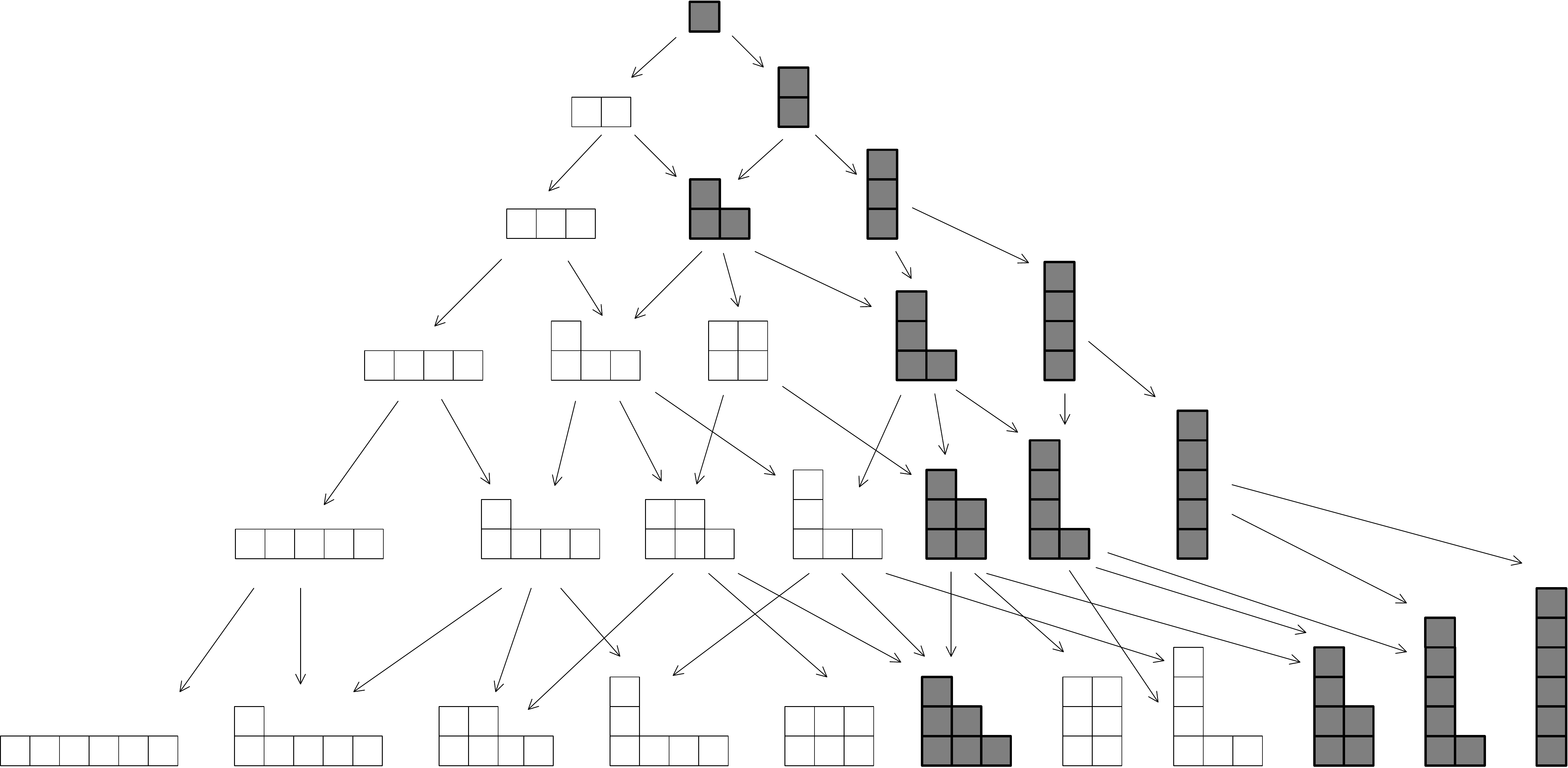}
	\caption{The beginning of the Young graph and the Schur graph.}
	\label{fig:schur}
	\end{figure}

A set of Young diagrams is a set of vertices of an infinite oriented graded graph known as \textit{Young graph}. Two Young diagrams are connected in this graph if one of them can be obtained from the other by adding a box. A subgraph of the Young graph which consists only of strict Young diagrams is called \textit{the Schur graph}. Figure \ref{fig:schur} shows the beginning part of the Young graph. The diagrams coloured dark grey form the beginning part of the Schur graph. \textit{The width of a diagram} is the number of its columns. \textit{The size of a diagram} is the number of its boxes. We consider a set of all strict diagrams of size $i$ as \textit{the $i$th level of the Schur graph}.

Let us consider the arrangement of integer numbers from $1$ to $n$ in diagram boxes, where $n$ is the diagram size. This arrangement is called \textit{Young tableau} if numbers grow by columns and by rows. It can be easily seen that the boxes are enumerated with respect to possible order of their additions. A Young diagram which corresponds to a Young tableau is called a \textit{shape} of the Young tableau. It is obvious that Young tableaux correspond to the paths from the root of the Young graph to the diagram of a given shape. The number of all paths to the diagram is called \textit{the dimension of a standard Young diagram}. \textit{The maximum diagram of a fixed size} is a diagram with the greatest possible dimension among all the diagrams of the corresponding size. In the same way, we define \textit{the dimension of a strict Young diagram} as a number of all paths to it from the root of the Schur graph.

\begin{figure}[h]
\captionsetup[subfigure]{labelformat=empty}
	\centering
	\subfloat[][]{\includegraphics[scale=0.6]{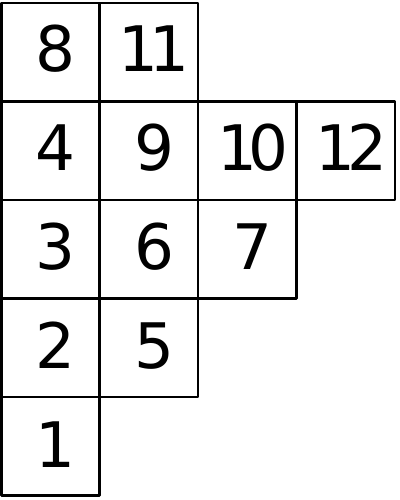}}
\hspace{100pt}
     	\subfloat[][]{\includegraphics[scale=0.6]{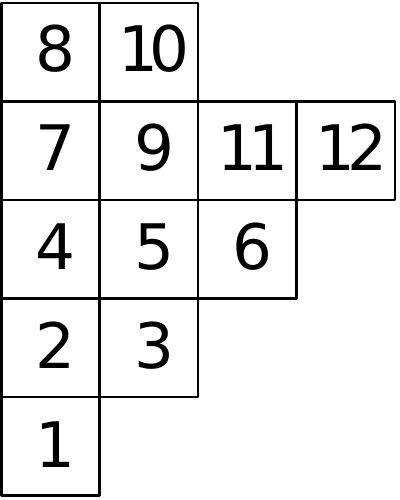}}
     	\caption{Examples of skew Young tableaux.}
     	\label{fig:skew}
	\end{figure}

Strict Young diagrams and the paths in the Schur graph are in one-to-one correspondence with so-called \textit{skew Young diagrams} and \textit{skew Young tableaux} respectively. Let us shift up each column of a strict Young diagram in such a way that the column number $i$ is shifted by $i-1$ boxes. A skew Young tableau is such an arrangement of integer numbers from $1$ to $n$, where $n$ is the diagram size, that the numbers grow by columns and by rows. It is obvious that such skew tableaux correspond to the paths in the Schur graph. Figure \ref{fig:skew} shows two examples of skew Young tableaux corresponding to the strict Young diagram from Fig.~\ref{fig:strict}.

The study of asymptotics of normalized dimensions is a subject of great interest in asymptotic representation theory. This asymptotics gives us information about the growth of weights in maximum irreducible representations of the symmetric group. The problem of investigation of this asymptotics has been formulated in \cite{verker85}.

The conjecture about normalized dimensions convergence for the case of typical Plancherel diagrams has also been raised in that paper. For the case of standard Young diagrams, this problem was studied in \cite{verpavl, ius, pdmi}. The results of experiments discussed in these papers give enough reason to consider that Vershik\,--\,Kerov conjecture is correct. Some rough estimations for the limits of normalized dimension can also be found there.

\subsection{Normalized dimensions of strict Young diagrams}

There are exact formulae to calculate the dimension of a diagram $\lambda$ for both Young graph and Schur graph cases. For Young graph, the dimension of a diagram $\lambda$ can be represented by the famous hook length formula:

\begin{equation}
\nonumber
\mathrm{dim}(\lambda) = \frac{n!}{\prod_{(i,j) \in \lambda} h(i,j)},
\end{equation}

where $n$ is the diagram size, $(i,j)$ are the coordinates of the box in the diagram, and $h(i,j)$ is the length of the corresponding hook, i.e. the set of boxes to the right and above the box, including the box itself.

The dimension of a strict Young diagram can be calculated in a more convenient way \cite{lpetrov, fpetrov, ivanov}:

	\begin{equation}
	\label{eq:strictdim}
	\mathrm{dim}(\lambda) = \prod_{i < j} \frac{l_i-l_j}{l_i+l_j} \cdot \frac{n!}{\prod{l_i!}},
	\end{equation}

where $l_i$ is the height of the $i$th column of the diagram and $n$ is the diagram size.
This formula can be derived from an analogous hook length formula for the number of skew Young tableaux by means of simple transformations. 

The dimensions may grow rapidly while the size of the diagrams increases. Therefore, it is crucial to use some appropriate normalization in order to study the asymptotic behaviour of the dimensions. Here we suggest the following normalization method:

	\begin{equation}
	\label{eq:normdim}
	c(\lambda) = -\frac{\ln{\mathrm{dim}(\lambda)} - \ln{\sqrt{n!}} + \frac{\ln 2}{2} \cdot n }{\sqrt{n}}.
	\end{equation}

We call $c(\lambda)$ \textit{a normalized dimension} of a diagram $\lambda$. It should be noted that when $n$ is fixed, the greater the normalized dimension, the smaller the exact dimension. Standard Young diagrams of size $n$ parameterize irreducible representations of the symmetric group $S_n$. Thereby, dimensions of diagrams are exactly the weights of the irreducible representations. Strict diagrams (or vertices of the Schur graph) of size $n$ parameterize projective representations of the group $S_n$. However, the correspondence is not straightforward in this case. A strict diagram $\lambda$ corresponds to one projective representation of the group $S_n$ when $\|\lambda \|-l(n)$ is even. On the contrary, if $\|\lambda \|-l(n)$ is odd, then a diagram $\lambda$ corresponds to two representations, whose dimensions differ by $1$. 

Projective representations of the symmetric group $S(n)$ correspond to linear representations of the so-called spin-symmetric group $T_n$ which is a central $\mathbb{Z}_2$ extension of the symmetric group $S(n)$:
$$ 1 \rightarrow \mathbb{Z}_2 \rightarrow T_n \rightarrow S_n \rightarrow 1.$$
Linear representations of the group $T_n$ are split into two parts. The first part is parameterized by the partitions of $n$ and corresponds to irreducible representations of the group $S_n$. The second part is parameterized by strict partitions. The detailed construction of these representations and their weights was discussed in \cite{nazarov}.

\subsection{Plancherel measure}

A central measure called Plancherel measure can be defined on the Schur graph as follows \cite{borodin}:

	\begin{equation}
	\nonumber
	Pl_n(\lambda) = \frac{\mathrm{dim}(\lambda)^2}{n!} \cdot 2^{n-l(\lambda)},
	\end{equation}

where $\mathrm{dim}(\lambda)$ is the dimension of the diagram $\lambda$, $n$ is the diagram size, $l(\lambda)$ is the width of the diagram.

This measure corresponds to a central Markov process named \textit{Plancherel process}. The property of centrality means that all the paths to the diagram $\lambda$ have the same probabilities. A path which appears under the Plancherel process is called \textit{Plancherel sequence}. 

The properties of paths or diagrams which have probability 1 and correspond to Plancherel measure are called \textit{Plancherel typical properties} of the paths or diagrams. Sometimes, we use the terms "typical diagram" or "typical sequence" in obvious sense.

\section{Greedy sequences and merging conjecture}

Let us consider a pair of Young diagrams $\lambda_1$, $\lambda_2$ such as $\lambda_2$ is obtained from $\lambda_1$ by adding a single box. Any path in the Young graph or the Schur graph consists of such elementary steps. An elementary step with a fixed $\lambda_1$ is called greedy if $\lambda_2$ has the maximum possible dimension. If each next diagram $\lambda_i, i>1,$ of a sequence $L=\lambda_1,\lambda_2,...,\lambda_N$ is obtained from $\lambda_{i-1}$ by a greedy step, we call $L$ \textit{a greedy sequence}. It is obvious that a greedy sequence can be started from any diagram $\lambda$.

Let us introduce so-called \textit{tail equivalency relations} on a set of infinite paths in the Young graph or the Schur graph. Two sequences of Young diagrams are tail-equivalent if they differ in a finite number of elements. Obviously, any asymptotic property of a sequence is a property of a whole equivalency class. 

There is an interesting property checked in many of our computer experiments. We have found that supposedly all greedy sequences belong to the same equivalency class. Briefly speaking, two greedy sequences always merge independently from their starting diagrams. It is true for greedy sequences in either Young graph \cite{ius} or Schur graph.

In order to verify this conjecture by computer experiments, we generated 5,000 pairs of typical Plancherel diagrams for 5,000 boxes. Greedy processes were started from each pair. All of them merged in a single sequence before 300 boxes were added. A similar single experiment was performed for two typical diagrams of 1 million boxes. The greedy processes merged after 23,466 steps.

\begin{figure}[!hb]
	\centering
		\subfloat[][]{\includegraphics[scale=0.77]{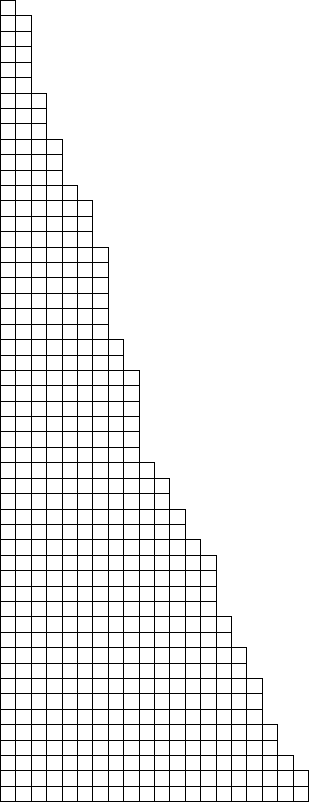}\label{fig:merging1}}
\hspace{2cm}
	     	\subfloat[][]{\includegraphics[scale=0.77]{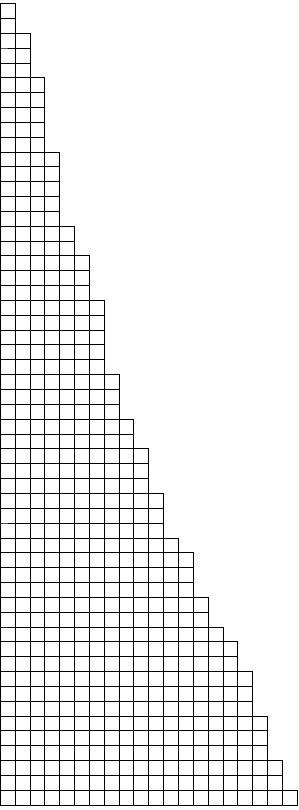}\label{fig:merging2}}
\hspace{2cm}
	     	\subfloat[][]{\includegraphics[scale=0.77]{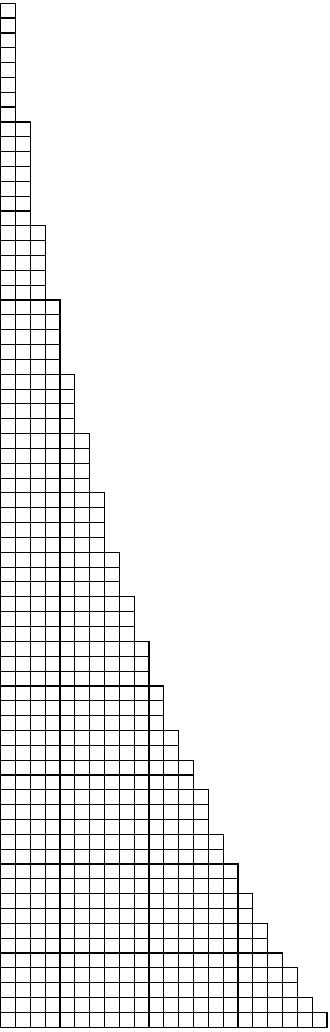}\label{fig:merging3}}
	     	\caption{(a),(b) Two initial diagrams of greedy sequences, (c) a diagram in a point where these sequences merge.}
	     	\label{fig:3diags}

	\end{figure}

Figures \ref{fig:merging1} and \ref{fig:merging2} show two Young diagrams of 500 boxes generated by random Plancherel process. A pair of greedy sequences started from these diagrams merges just after 69 steps. Their joint diagram of 569 boxes in a merge point is depicted in Fig.~\ref{fig:merging3}.

In the case when the initial diagrams significantly differ in their shapes, the number of steps required for the merging is much greater. We studied such extreme cases when the initial diagrams differ a lot. We chose a single column diagram and a triangle diagram (where the heights of columns differ by one, and the rightmost column consists of one box) as initial diagrams. We investigated such extreme pairs of diagrams consisting up to 518 boxes. In all our computer experiments, the greedy sequences which started from these pairs merged. For example, the greedy sequences which started from the single column diagram and the triangle diagram consisting of 518 boxes merged after 34,697 steps.

This observation gives us a reason to hypothesize that all greedy sequences of strict Young diagrams contain infinite number of maximum dimension diagrams. It can be easily seen that if the previous statement is true, then the asymptotics of diagram dimensions in a greedy sequence is equal to the asymptotics of maximum dimensions. On the contrary, typical sequences built by Plancherel process often have no joint diagrams. 

Originally, our motivation to study greedy sequences was that any greedy sequence allows us to investigate the exact asymptotics of maximum dimensions. As it turned out, greedy sequences have many other interesting properties as well. Some of them are discussed in Sec.~\ref{sec:myst}. This merging conjecture was the very last topic discussed by the authors with Sergey V. Duzhin.

\section{Maximum diagrams searching algorithm}
\label{sec:algorithm}

Let us introduce the following order relation on a set of sequences of Young diagrams. We say that a sequence of diagrams $X=(x_1,\dots,x_{n})$ where $\mathrm{size}(x_i) = \mathrm{size}(x_{i-1}) + 1$ is bigger than a similar sequence $Y=(y_1,\dots,y_{m})$ if for a certain $k\leq \min(n,m)$ we have $\dim x_i=\dim y_i$ for all $i<k$ and $\dim x_k>\dim y_k$. Note that generally $x_{i+1}$ is not necessary obtained from $x_{i}$ by adding one box.

It can be easily seen that under such a definition, the sequence of diagrams with maximum dimensions is maximum among all possible sequences. We say that a bigger sequence is better than a smaller one. Briefly speaking, the bigger sequence improves the smaller one. Earlier in \cite{pdmi} an algorithm which improves any given sequence was introduced. In this paper, we propose a much more efficient technique for constructing sequences of diagrams with large dimensions. 

The main idea of the algorithm is the following. At a certain level of the Schur graph, we choose some set of diagrams and start greedy sequences from them. At each next level, we remember the diagram with the highest dimension among all the diagrams of the corresponding size in all the sequences under the process. It is obvious that if one of these sequences passes through some maximum diagram, this diagram will be remembered in the resulting sequence. Note that these several sequences merge into a single one after a small number of steps. When these sequences are merged, we repeat the previous steps starting from a higher level of the Schur graph. Since we start greedy sequences from different diagrams, we obtain different sequences which may have strictly maximum diagrams at different levels. The resulting sequence contains all of these maximum diagrams. We choose the diagram with maximum dimension in each path, and since the algorithm is fed by the best constructed diagrams, we hope it finds many exact maximum diagrams. This kind of process consistently improves any sequence of diagrams which contains non-maximum diagrams.

We executed this algorithm several times and constructed a sequence of diagrams with large dimensions up to 500,000 boxes. The first diagrams of this sequence match all known exact maximum diagrams. None of our previous algorithms \cite{pdmi} can improve any diagrams of this sequence. We suppose that most of the diagrams in this sequence are strictly maximum. Figure \ref{fig:greedy_and_improved} illustrates the comparison of the normalized dimensions of diagrams from a greedy sequence and from the sequence made by this algorithm. Nevertheless, it can be seen from this figure that some of normalized dimensions in the greedy sequence are strictly equal to the ones in the improved sequence.

\section{Study of normalized dimension asymptotic behaviour in greedy sequences}

\begin{figure}[!htb]
	\centering
	\includegraphics[scale=0.4]{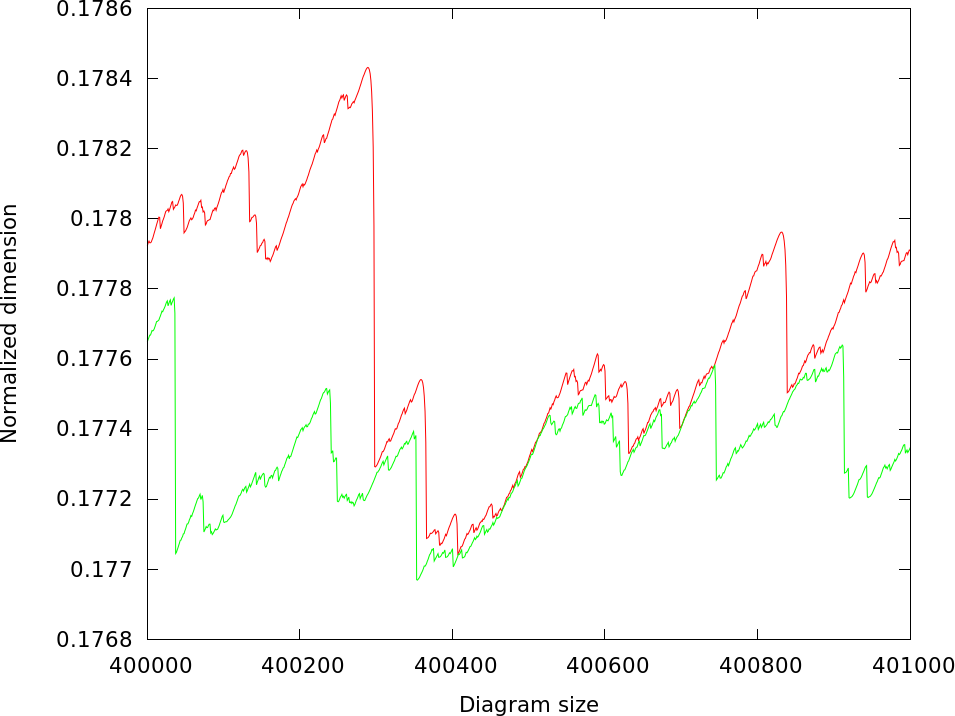}
	\caption{Normalized dimensions of strict diagrams in a greedy sequence (\protect\includegraphics{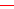}) and in an improved sequence (\protect\includegraphics{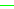}).}
	\label{fig:greedy_and_improved}
	\end{figure}

In contrast with the previous papers \cite{ius, pdmi}, our study of asymptotic behaviour of normalized dimensions as function of the diagram size is based on the investigation of their finite differences. We study not only the growth but also the oscillations of normalized dimensions in greedy sequences. It can be done by studying finite differences or a discrete derivative of a sequence of normalized dimensions. Figure \ref{fig:greedy_and_improved} shows that normalized dimensions are represented by a fast-oscillating function.
Sometimes, we need to work with exact dimensions, for example, when we compare two diagrams with normalized dimensions close to each other. But these exact dimensions are so large that operations with them are extremely complex from the computational point of view. To do that, we compute the corresponding sequence of exact dimension ratios for each sequence of diagrams. This trick allows us to compute normalized dimensions with great accuracy and to compare exact dimensions.

\section{The finite differences of normalized dimensions in greedy sequences \,--\, the mysterious properties.}
\label{sec:myst}

\subsection{The case of strict diagrams}

\begin{figure}[!h]
	\centering
	\subfloat[][]{\includegraphics[width=.5\textwidth]{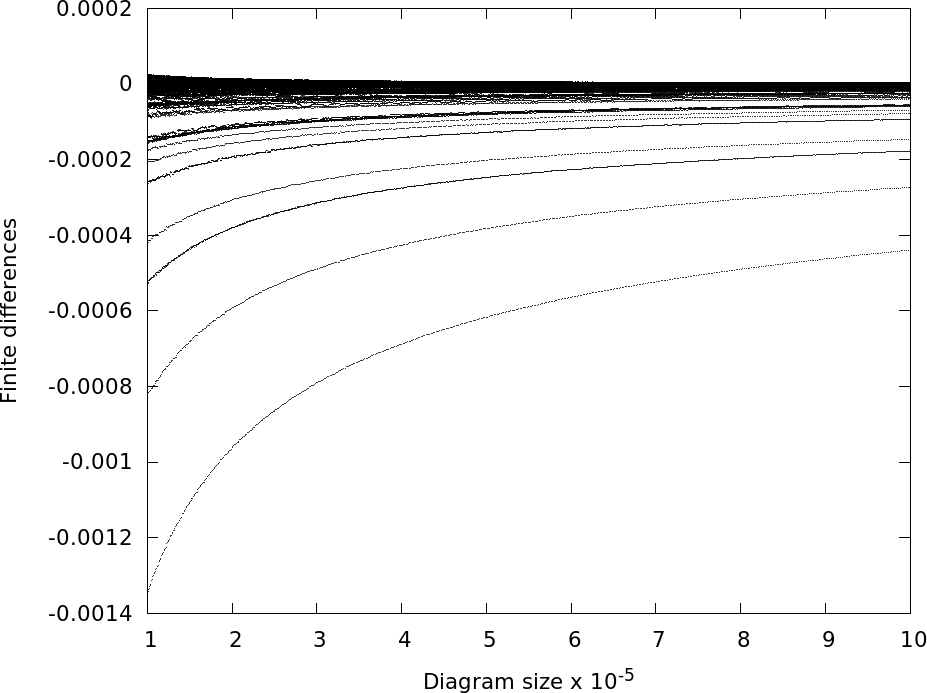}\label{fig:strict_greedy_diffs}}
\hspace*{\fill}
     	\subfloat[][]{\includegraphics[width=.5\textwidth]{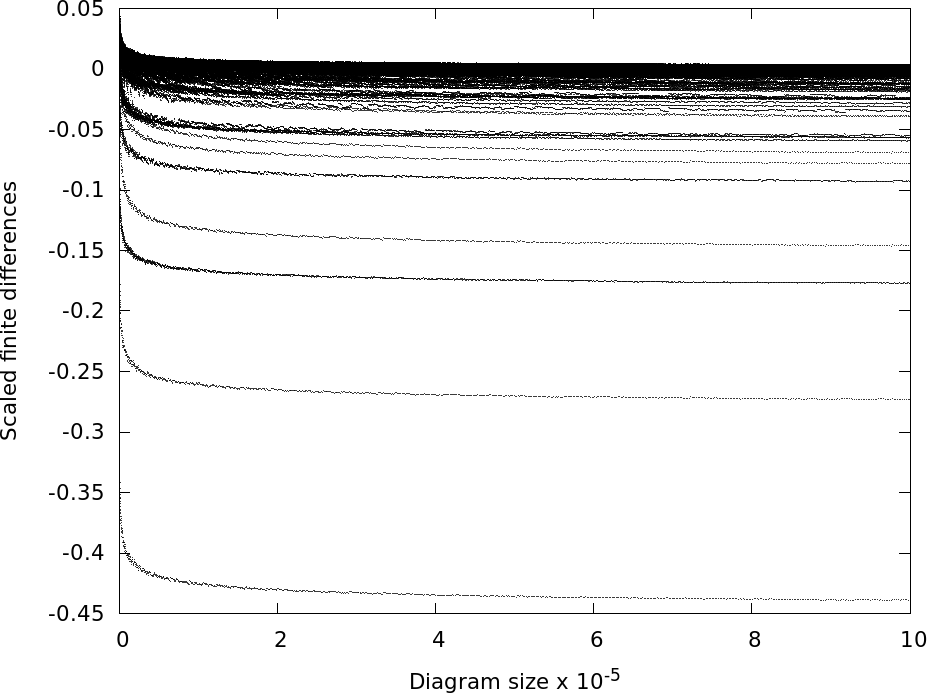}\label{fig:strict_greedy_diffs_sqrtn}}
     	\caption{Finite differences of normalized dimensions in a greedy sequence of strict diagrams: (a) original function, (b) multiplied by $\sqrt{n}$.}
	\end{figure}

Let us consider a sequence of differences between normalized dimensions $c(n)-c(n-1)$ of neighbouring diagrams in greedy sequence of strict Young diagrams. The next diagram is obtained from the previous one by adding one box. The graph of these differences is shown in Fig.~\ref{fig:strict_greedy_diffs}. Note that this picture represents the usual graph of the function that contains $10^6$ points. For each point, we know the coordinates of the corresponding box.

Analysis of this graph allows us to understand the structure of oscillation of normalized dimensions. The positive values in Fig.~\ref{fig:strict_greedy_diffs} correspond to the rises of the normalized dimension function, while the negative values correspond to the falls. The maximum falls are represented on the graph as lines.

Based on behaviour analysis of these finite differences from Fig.~\ref{fig:strict_greedy_diffs}, we suggest a conjecture that the decay rate of the oscillations is of order $1/\sqrt{n}$. We multiply these finite differences by $\sqrt{n}$ in order to investigate the sequence $(c(n)-c(n-1)) \cdot \sqrt{n}$. Further we say "differences" instead of "differences $\cdot \sqrt{n}$" if the meaning is clear from the context. The sequence of these differences is depicted in Fig.~\ref{fig:strict_greedy_diffs_sqrtn}. We discuss some interesting observations of this sequence below.

First, it can be seen that all these differences are limited by range [-0.44, 0.3]. The interval of values for large $n$ is even narrower: for $n>2 \cdot 10^5$ it is [-0.44, 0.01]. The differences at $n = 2\cdot10^5$ and $n = 10^6$ do not differ much. Hence, we suppose that this behaviour does not change for $n \rightarrow \infty$.

\begin{figure}[!htb]
	\centering
	\includegraphics[scale=0.4]{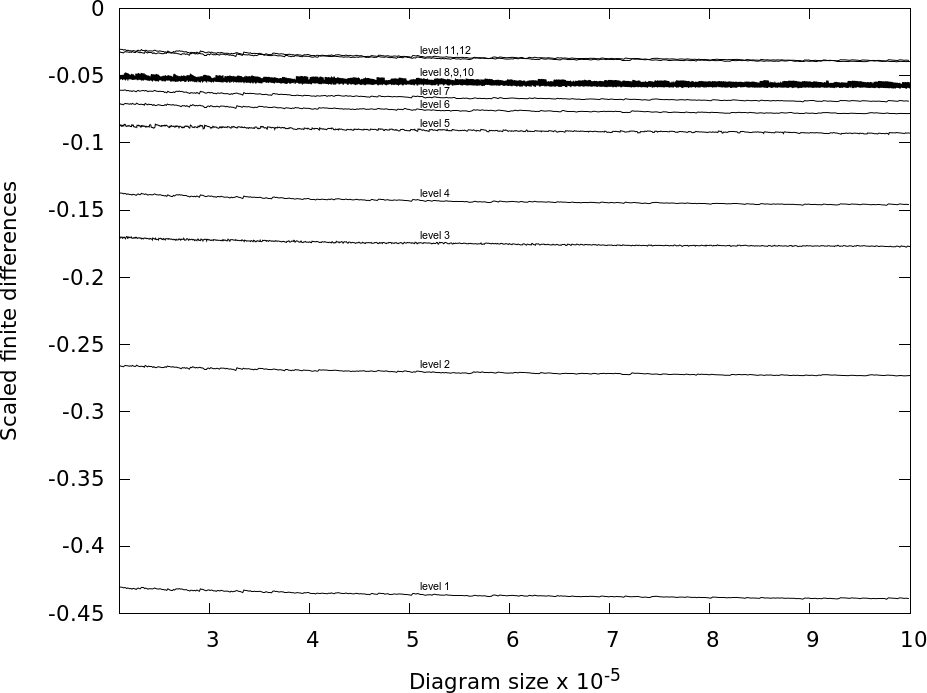}
	\caption{First 12 levels.}
	\label{fig:first_12_levels}
	\end{figure}

Next, we can see from this picture that the values of the normalized dimension differences are located as follows. The density of these values is high in the area close to zero. At the same time, we can observe a certain number of discrete levels represented by some negative values located in a set of narrow intervals. It is important to note that there are no values in the spaces between the level intervals. This picture is mysterious in many ways, because the meanings of different levels are not the same. The bottom 12 levels are depicted separately in Fig.~\ref{fig:first_12_levels}. Let us enumerate the levels upwardly: the number of the lowest level will be 1 and so on.

It was found that all values of the first level starting from $n > 10200$ are in one-to-one correspondence to the diagrams with the last box added by a greedy step in the 7th row. It means that in the greedy sequence of strict diagrams, the largest falls of normalized dimension happen when a box is added to the row number 7. This fact looks very strange. Similarly, the values of the next level match an addition to the row number 2. Table~\ref{table:levels_rows} shows some other similar examples of associations between the level number and the number of the row. 

\begin{table}[ht]
\centering
	\caption{The correspondence of levels to the numbers of rows to which boxes were added.\vspace{.5em}}
	   \label{table:levels_rows}
		\begin{tabular}{@{}|c|c|c|c|c|c|c|c|c|c|c|c|c|}
		\hline
			No. of level & 1 & 2 & 3      & 4 & 5      & 6 & 7 & 8      & 9      & 10     & 11 & 12\\ 
		\hline
			No. of row   & 7 & 2 & $\ast$ & 5 & $\ast$ & 4 & 3 & $\ast$ & $\ast$ & $\ast$ & 6  & 9 \\
		\hline
		\end{tabular}
\end{table}

\begin{figure}[!htb]
   \begin{minipage}{0.49\textwidth}
     \centering
     \includegraphics[width=\textwidth]{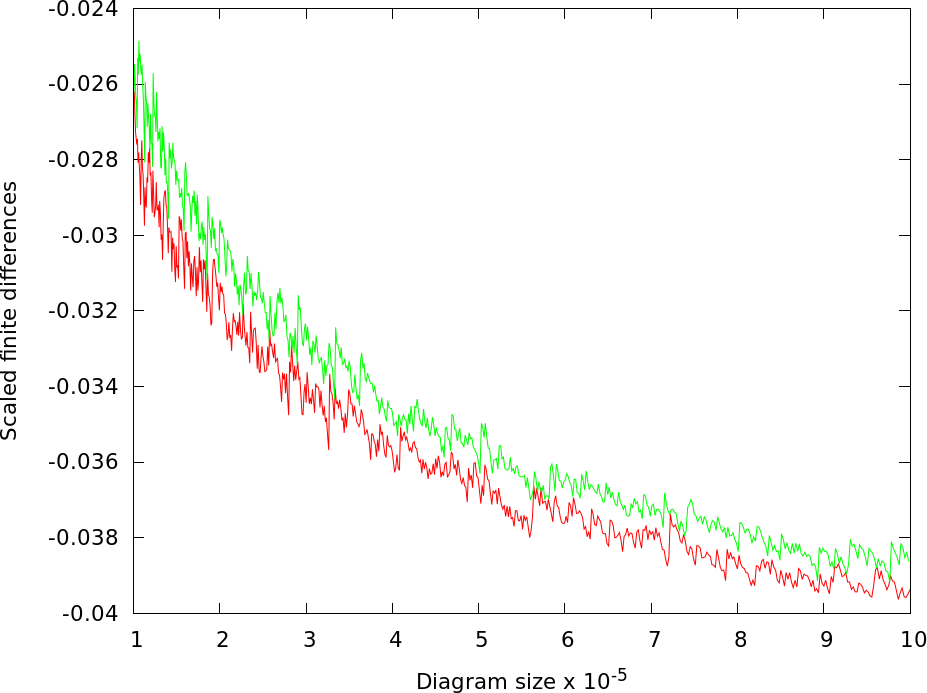}
	\caption{The differences of normalized dimensions multiplied by $\sqrt{n}$ when a box is added to the 6th row (\protect\includegraphics{line_red.png}) and the 9th row (\protect\includegraphics{line_green.png}).}
	\label{fig:to_6_and_9}
   \end{minipage}\hfill
   \begin{minipage}{0.49\textwidth}
     \centering
	\includegraphics[width=\textwidth]{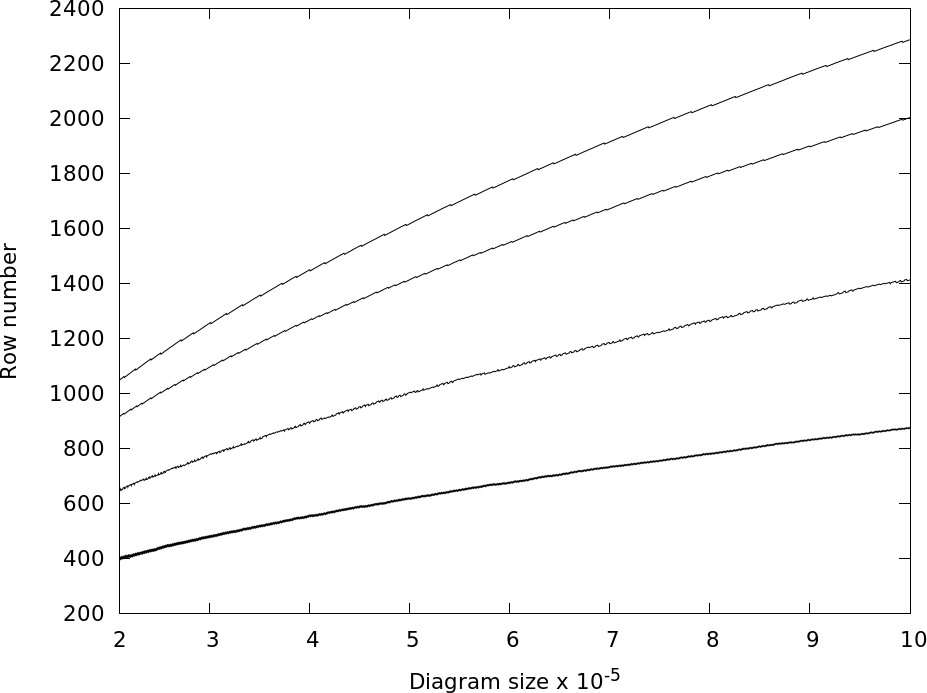}
	\caption{Some special levels and the corresponding row numbers.\\~}
	\label{fig:curves}
  \end{minipage}
\end{figure}

The values for the 11th and 12th levels are very close to each other, so they are shown in Fig.~\ref{fig:to_6_and_9} in a larger scale.
However, not every level corresponds to a fixed number of a row where a box was added in a greedy algorithm. There are also levels of a different type which are denoted by stars in Table~\ref{table:levels_rows}. For example, the values from the interval [-0.178, -0.175], represented at the level 3 (see Fig.~\ref{fig:first_12_levels}), correspond to the boxes added to rows with growing numbers. The consecutive points located at this level correspond to the sequence of added boxes with ascending numbers of both rows and columns. Figure \ref{fig:curves} shows how an added row number depends on the size of a diagram for some levels of this type.
We expect that some part of the differences located close to zero consists of a large set of curves which are similar to those presented in Fig.~\ref{fig:curves}.

\begin{figure}[!htb]
	\centering
	\includegraphics[scale=0.4]{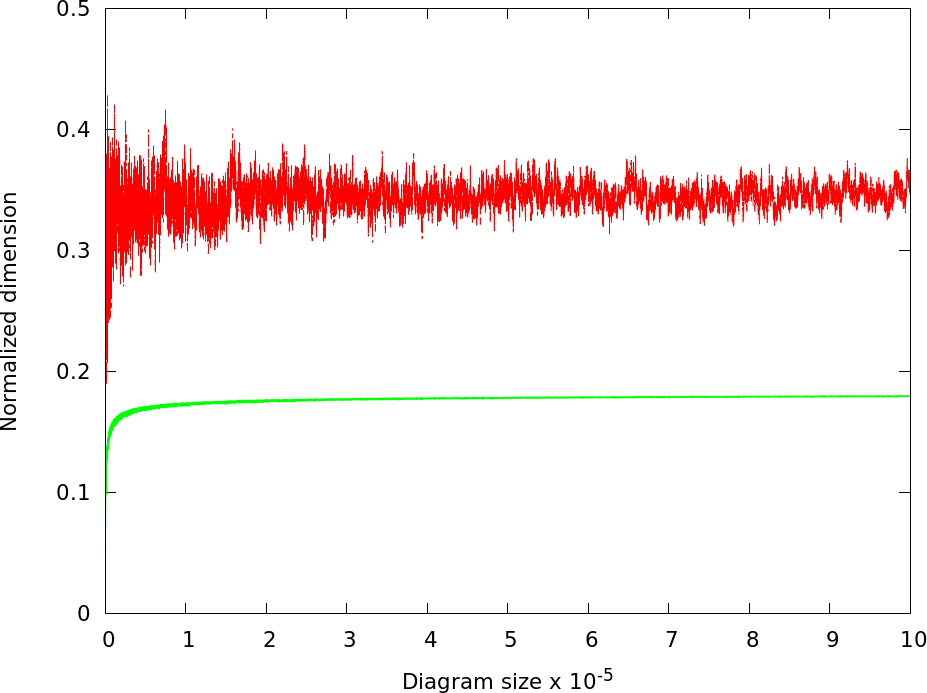}
	\caption{Normalized dimension of strict diagrams in the greedy sequence (\protect\includegraphics{line_green.png}) and in a typical Plancherel sequence (\protect\includegraphics{line_red.png}).}
	\label{fig:greedy_planch}
	\end{figure}

It is known \cite{verker85} that maximum dimension diagrams and typical Plancherel diagrams have the same limit shape. But normalized dimensions in these cases differ from each other significantly (see Fig.~\ref{fig:greedy_planch}). 
So, it is necessary to understand which geometrical property makes these sequences different. In order to study this phenomenon, we investigated the partitions which correspond to exact maximum Young diagrams of size up to 250 boxes obtained by exhaustive search. Besides, we investigated the partitions corresponding to an improved sequence constructed by the proposed algorithm.

Let us call the rightmost part of a partition \textit{a regular tail} if it contains only consecutive odd (..., 7, 5, 3, 1) or consecutive even (..., 8, 6, 4, 2) numbers. We have found out that most of partitions corresponding to maximum diagrams have a regular tail of a relatively large size. Some examples of such partitions are listed in Table~\ref{table:max_partitions}. 

\begin{table}[ht]
	\centering
	\caption{Examples of partitions for some maximum diagrams.}
	    \label{table:max_partitions}
   	 \begin{tabular}{|@{}c|c@{}|}
		\hline
	    Size & Partition  \\
		\hline
	    200 & \hspace{1em} 34, 30, 26, 23, 20, 17, 14, \textbf{11, 9, 7, 5, 3, 1}\\
		\hline
	    250 & 38, 34, 30, 27, 24, 21, 18, 15, 13, \textbf{10, 8, 6, 4, 2}\\
		\hline
	    \end{tabular}
\end{table}

The regular tail is shown in bold. We expect that this property holds for all sufficiently large maximum diagrams. This property also holds for most strict Young diagrams from the sequence constructed by our algorithm described in Sec.~\ref{sec:algorithm}. For example, the regular tail of the partition corresponding to the diagram of size = 100,000 is (47, 45, 43, 41, 39, 37, 35, 33, 31, 29, 27, 25, 23, 21, 19, 17, 15, 13, 11, 9, 7, 5, 3, 1). If the diagram has a size equal to 300,000, the regular tail of its partition will be (64, 62, 60, 58, 56, 54, 52, 50, 48, 46, 44, 42, 40, 38, 36, 34, 32, 30, 28, 26, 24, 22, 20, 18, 16, 14, 12, 10, 8, 6, 4, 2). 

\subsection{The case of standard diagrams}
	\begin{figure}[!htb]
	\centering
	\includegraphics[scale=0.4]{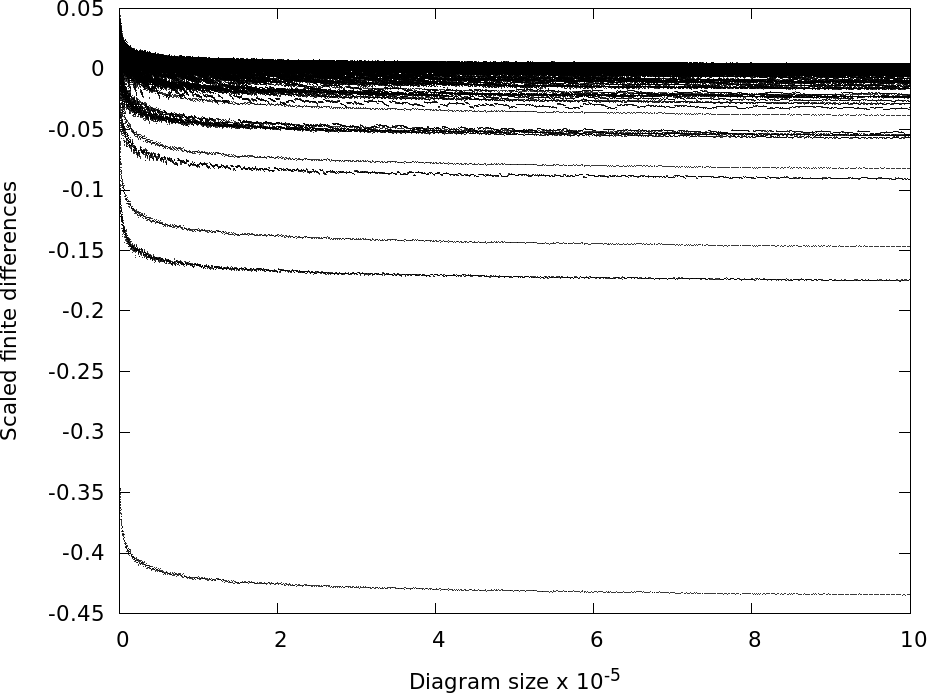}
	\caption{Finite differences of normalized dimensions in a greedy sequence of standard diagrams multiplied by $\sqrt{n}$.}
	\label{fig:standard_diffs_sqrtn}
	\end{figure}

A similar analysis of the finite difference function can also be made for the case of standard diagrams. In this case, the formula for normalized dimension of diagram $\lambda$ is as follows:
$$
c(\lambda) = -\frac{1}{\sqrt{n}} \ln \frac{\mathrm{dim} \lambda}{\sqrt{n!}}.
$$
Note that this notation differs from the notation used in \cite{verpavl, ius}.

The distribution of values $(c(n)-c(n-1)) \cdot \sqrt{n}$ in a greedy sequence of standard diagrams is shown in Fig.~\ref{fig:standard_diffs_sqrtn}. It looks very similar to the case of strict diagrams (see Fig.~\ref{fig:strict_greedy_diffs_sqrtn}). However, each of the bottom levels corresponds to the differences of the added row and column numbers instead of just a row number. We will not go into details about the standard diagram case because it is similar to the case of strict diagrams discussed above.

\subsection{Estimations of limits for normalized dimensions}

The fact that the differences of normalized dimensions multiplied by $\sqrt{n}$ are located within the interval [-0.44, 0.3] in the case of a greedy sequence of strict Young diagrams makes it possible to estimate the degree of convergence of normalized dimensions to some limit value.

As the differences of normalized dimensions multiplied by $\sqrt{n}$ are limited, we make an assumption that for a hypothesized limit value $C$ the function of differences $(C-c(n)) \cdot \sqrt{n}$ is limited, too. This assumption allows us to give a much more precise estimation of the normalized dimension limit.

The idea about this adjustment is finding the best mean-square approximation of the function $c(n)$ with a function of the form $C + a/\sqrt{n}$ on a sufficiently large interval, wherein the computed value $C$ gives us an estimation for the limit of the function $c(n)$. This conjecture must imply that the degree of convergence of $c(n)$ to $C$ is of the order $\frac{O(1)}{\sqrt{n}}$. Actually, this difference $C-c(n)$ is close to $-\frac{4.29412}{\sqrt{n}}$ for the standard diagram case and to $-\frac{3.19579}{\sqrt{n}}$ for the strict diagram case.

This estimation method is much more precise than the previously used ones. It works for both standard and strict Young diagrams. 
The estimations obtained by this technique are as follows: 0.700281 for standard diagrams and 0.182837 for strict diagrams.
Note that the very existence of the limit still remains unproven.

\begin{figure}[!h]
	\centering
	\subfloat[][]{\includegraphics[width=.5\textwidth]{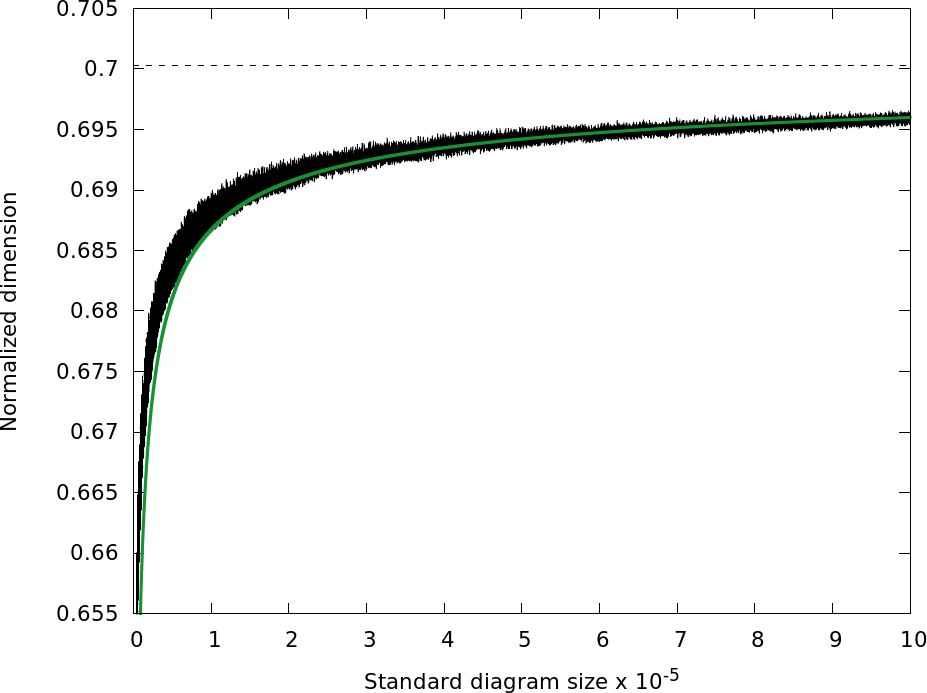}\label{fig:standard_lim}}
\hspace*{\fill}
     	\subfloat[][]{\includegraphics[width=.5\textwidth]{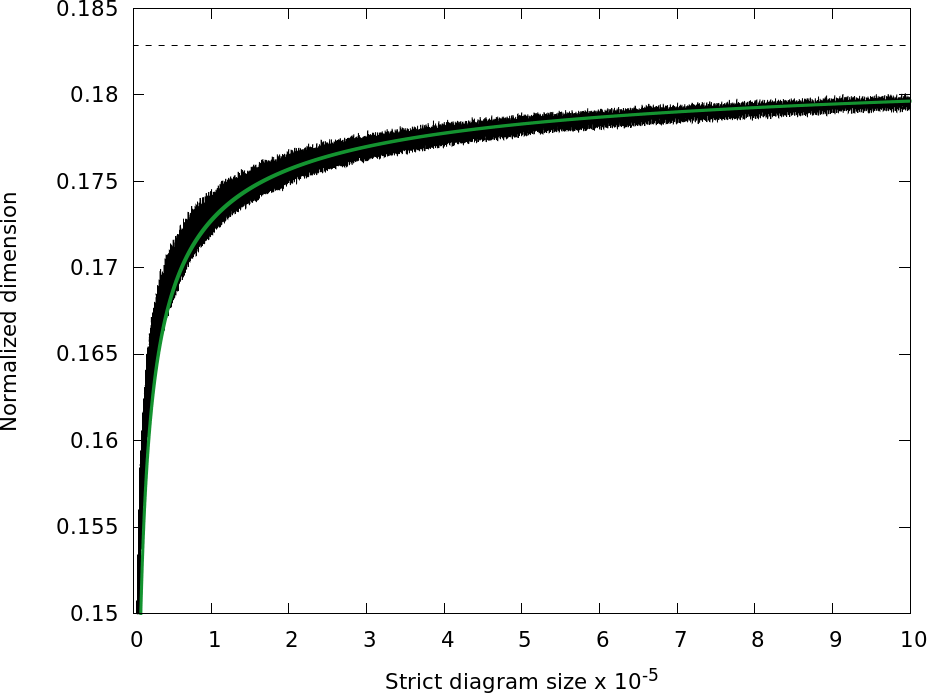}\label{fig:strict_lim}}
     	\caption{Convergence of normalized dimensions to a possible limit for (a) the greedy sequence of standard diagrams, (b) the greedy sequence of strict diagrams.}
     	\label{fig:limits}
	\end{figure}

Figure \ref{fig:standard_lim} shows the limit of normalized dimensions (dashed line), the graph of normalized dimensions (black solid line) and the function $0.700281-\frac{4.29412}{\sqrt{n}}$ (green solid line).

It can be seen that these graphs agree very well with our conjecture of degree of convergence.


\section*{Acknowledgments}
This work was supported by grant RSF 14-11-00581. The authors thank Vadim Smolensky for his comments that greatly improved the manuscript.

\end{document}